 0
 0

\def\Per{{\mathcal P}}
\def\F{{\mathcal F}}
\def\FF{{\mathfrak F}}
\newcommand{\fss} {$\FF$-semisimple}

\newcommand{\Aut} [1] { {\rm{Aut(}}#1 {)} }

\newcommand{\Ad} [1] { {\rm Ad}_{#1} }

\newcommand{\R} [1] { {\bf{R}}^{#1} }
\newcommand{\T} [1] { {\bf{T}}^{#1} }
\newcommand{\Z} [1] { {\bf{Z}}^{#1} }

\newcommand{\Q} [1] { {\bf{Q}}^{#1} }
\newcommand{\C} [1] { {\bf{C}}^{#1} }

\newcommand{\B} [0] { {\mathbb{B}} }

\newcommand \im { {\rm im}\, }

\newcommand {\E} [0] { {\bf E} }
\newcommand {\G} [0] { {\bf G} }

\newcommand {\A} [0] { {\mathfrak a} }

\documentclass{elsart}
\usepackage{amssymb}

\begin{document}

\begin{frontmatter}
\title{Invariant fibrations of geodesic flows \\
{\tt \,  pi1.tex}}

\author[LTB]{Leo T. Butler}
\address{Department of Mathematics and Statistics\\
Queen's University,
Kingston, ON, Canada
K7L 3N6}
\ead{lbutler@mast.queensu.ca}
\thanks[LTB]{Thanks to Keith Burns, Oleg Bogoyavlenskij and a referee
for their comments on this paper. Research partially supported by the
Natural Sciences and Engineering Research Council of Canada.}

\begin{abstract} 
Let $(\Sigma,{\bf g})$ be a compact $C^2$ finslerian $3$-manifold. If
the geodesic flow of ${\bf g}$ is completely integrable, and the
singular set is a tamely-embedded polyhedron, then $\pi_1(\Sigma)$ is
almost polycyclic. On the other hand, if $\Sigma$ is a compact,
irreducible $3$-manifold and $\pi_1(\Sigma)$ is infinite polycyclic
while $\pi_2(\Sigma)$ is trivial, then $\Sigma$ admits an analytic
riemannian metric whose geodesic flow is completely integrable and
singular set is a real-analytic variety.  Additional results in higher
dimensions are proven.
\end{abstract}


\end{frontmatter}

\section{Introduction}

A smooth ($C^1$) flow $\phi_t : M \to M$ is {\em integrable} if there
is an open, dense subset $R$ that is covered by angle-action charts
$(\theta,I) : U \to \T{k} \times \R{l}$ which conjugate $\phi_t$ with
a translation-type flow $(\theta,I) \mapsto (\theta + t \omega(I),
I)$. Evidently, there is an open dense subset $L \subset R$ fibred by
$\phi_t$-invariant tori~\cite{Bogo:1997a}. Let $f : L \to B$ be the
$C^1$ fibration which quotients $L$ by these invariant tori and let
$\Gamma = M - L$ be the {\em singular set}. If $\Gamma$ is a
tamely-embedded polyhedron, then $\phi_t$ is called $k$-{\em
semisimple} with respect to $(f,L,B)$. We say $\phi_t$ is semisimple
if it is $k$-semisimple with respect to some $(f,L,B)$.

A geodesic flow on a unit-sphere bundle $S \Sigma$ is {\em completely
integrable} if it is integrable with invariant tori of dimension $\dim
\Sigma$. It is evident that $\dim \Sigma$-semisimplicity is a
definition of topologically-tame complete integrability.

A group $\Delta$ is {\em polycyclic} of step length $c$ if there is a
finite chain of subgroups $1 = \Delta_0 \lhd \cdots \lhd \Delta_c =
\Delta$, with $\Delta_{i-1}$ normal in $\Delta_i$ and
$\Delta_i/\Delta_{i-1}$ cyclic for all $i$.  In the present paper,
$\Sigma$ is always a boundaryless manifold.

\begin{thm} \label{thm:1}
Let $(\Sigma,{\bf g})$ be a compact $C^2$ finslerian $3$-manifold. If
the geodesic flow of ${\bf g}$ is $3$-semisimple, then $\pi_1(\Sigma)$
contains a finite-index polycyclic subgroup of step length at most
$4$.
\end{thm}

Recall that $\Sigma$ is {\em irreducible} if every tamely-embedded
$2$-sphere bounds a $3$-ball. Evans and Moser classify the solvable
groups that appear as the fundamental group of a compact $3$-manifold
in~\cite{EvansMoser:1972a}; in this case they are all
polycyclic. Their result, along
with~\cite{BolsinovTaimanov:2000a,Butler:1999a}, help to prove

\begin{thm} \label{thm:2}
Let $\Sigma$ be a compact, irreducible $3$-manifold with
$\pi_1(\Sigma)$ infinite polycyclic and $\pi_2(\Sigma) = 0$. Then
$\Sigma$ admits an analytic riemannian metric whose geodesic flow is
$3$-semisimple.
\end{thm}

Four comments on Theorem~\ref{thm:2}: First, irreducibility is a
technical hypothesis that precludes $\Sigma$ from containing fake
$3$-balls. The Poincar\'e conjecture would make redundant this
hypothesis. Second, Evans and Moser's work, along with that of Hempel
and Jaco~\cite{HempelJaco}, implies that if $\Sigma$ satisfies the
hypotheses of Theorem~\ref{thm:2}, then $\Sigma$ admits flat geometry,
$Nil$-geometry or $Sol$-geometry. These geometries supply the
mentioned metrics. Third, if $\pi_1(\Sigma)$ is finite, then the
geometrization conjecture implies that $\Sigma$ admits
$S^3$-geometry. Fourth, if $\pi_1(\Sigma)$ is infinite polycyclic and
$\pi_2(\Sigma) \neq 0$, then Theorem 5.1 of~\cite{EvansMoser:1972a}
implies that $\pi_1(\Sigma)$ is isomorphic to one of $\pi_1(S^2 \times
S^1)$, $\pi_1(P^2 \times S^1)$ or $\pi_1(P^3 \# P^3)$. The universal
covering space $\tilde{\Sigma}$ of $\Sigma$ is two-ended and,
according to Ian Agol, the geometrization conjecture implies that
$\tilde{\Sigma}$ must be homeomorphic to $S^2
\times \R{}$. Some work establishes that $\Sigma$ itself is a
geometric $3$-manifold homeomorphic to $P^2 \times S^1$, $P^3
\# P^3$ or one of the two $S^2$-bundles over $S^1$. From the proof of
Theorem~\ref{thm:2} follows

\begin{thm} \label{thm:2.1}
Assume the geometrization conjecture. Let $\Sigma$ be a compact
$3$-manifold such that $\pi_1(\Sigma)$ is polycyclic. Then $\Sigma$
admits an analytic riemannian metric whose geodesic flow is
$3$-semisimple.
\end{thm}

Say $(f',L',B')$ is a {\em refinement} of $(f,L,B)$ if $B'$ is an open
dense subset of $B$, $L' = f^{-1}(B')$, $f' = f|L'$ and $\Gamma' = M -
L'$ is a tamely-embedded polyhedron; it is {\em tractable} if,
for each component $L'_i$ of $L'$, either (1) there is an
$f$-saturated, codimension-1 submanifold $W_i \subset L'_i$ such that
the inclusion map $\iota_{W_i,L'_i}$ is epimorphic on $\pi_1$ or (2)
there is a component $L'_j$ satisfying (1) and a map $r : L'_i \to
L'_j$ such that $\iota_{L'_i,M}$ is homotopic to $\iota_{L'_j, M}
\circ r$. By Lemma~\ref{lem:niceref}, $(f,L,B)$ has a tractable
refinement, so it can be assumed from the outset that $(f,L,B)$ is
tractable.

\noindent
A group is {\em small} if it does not contain a free group
on two generators.

\begin{thm} \label{thm:3}
Let $(\Sigma,{\bf g})$ be a compact $C^2$ finslerian $4$-manifold. If
the geodesic flow of ${\bf g}$ is $4$-semisimple and the fundamental
group of each component of $B$ is small, then $\pi_1(\Sigma)$ contains
a finite-index polycyclic subgroup of step length at most $6$.
\end{thm}

Theorems~\ref{thm:1} and~\ref{thm:3} have similar proofs. Since
$(f,L,B)$ is tractable, each component $L_i$ of $L$ ``looks like'' a
$\T{k}$-bundle over a codimension-1 submanifold $N \subset B$. $N$ is
a $1$-manifold in Theorem~\ref{thm:1} and a $2$-manifold whose
fundamental group is small in Theorem~\ref{thm:3}. In both cases,
$\pi_1(N)$ is almost abelian, which implies $\pi_1(L_i)$ is almost
polycyclic. Theorems~\ref{thm:1} and~\ref{thm:3} then follow from
Lemma~\ref{lem:2}, which states that there is a component $L_i$ such
that the inclusion map $\iota_{L_i, S\Sigma}$ has a finite-index image
in $\pi_1(S\Sigma)$. Similar ideas appear in Ta\u{\i}manov's
work~\cite{Taimanov:1988a}.

Without the condition on $B$ in Theorem~\ref{thm:3}, it seems
difficult to deduce any constraint on $\pi_1(\Sigma)$ with the
techniques of this paper. One is led to ask

\medskip

\noindent
{\bf Question A}: {\em is there any obstruction to the existence of a
$4$-semisimple geodesic flow on a $4$-manifold?}

\noindent
and more pointedly,

\noindent
{\bf Question B}: {\em is there a $4$-semisimple geodesic flow on $(S^3 \times
S^1) \# (S^3 \times S^1)$?}

\medskip

\noindent
{\em Algebraic properties of $\pi_1$ and recurrence}: A subgroup is
{\em almost normal} if its normalizer is of finite-index.

\begin{defn} \label{def:1}
A group is {\em anabelian} if its only abelian, almost
normal subgroup is the trivial group $1$.
\end{defn}

A Gromov-hyperbolic group is either anabelian or almost cyclic. Also,
the fundamental group of a finite-volume riemannian manifold of
non-positive curvature is either anabelian or almost abelian.

\def\Lim{{\mathcal L}}

Recall two notions of recurrence: (1) $x$ is periodic if there is a
$T>0$ such that $\phi_T(x) = x$. $\Per(\phi)$ is the closure of the
set of periodic points. (2) the $\omega$-limit set of $x$,
$\omega(x),$ is the set of points $y$ such that $\phi_{t_k}(x) \to y$
for some sequence $t_k \to +\infty$. The $\alpha$-limit set,
$\alpha(x)$, is defined similarly. If $x \in \omega(x) \cap
\alpha(x)$, then $x$ is a recurrent point. The limit-point set
$\Lim(\phi)$ is the closure of the recurrent-point set.

Let $\tilde{\Sigma}$ be the universal cover of $\Sigma$. The pullback
of $\bullet$ on $S\Sigma$ to $S\tilde{\Sigma}$ is $\bar{ \bullet }$.

\begin{thm}  \label{thm:4}
Let $(\Sigma,{\bf g})$ be a compact $C^2$ finsler manifold whose
geodesic flow $\phi_t$ is semisimple. If $\pi_1(\Sigma)$ is anabelian,
then ${\rm Int}\, \Lim(\bar{\phi}) \neq
\emptyset$. \\ In any $C^2$ neighbourhood of ${\bf g}$,
there is a finsler metric ${\bf k}$ with geodesic flow $\varphi_t$
such that ${\rm Int}\, \Per(\bar{\varphi}) \neq \emptyset$, and
$\phi_t$ {\rm \&} $\varphi_t$ share the same invariant tori.
\end{thm}

Following~\cite{Kozlov:1996a}, one can construct a $C^{\infty}$
riemannian metric ${\bf g}$ on $\Sigma$ with ${\rm Int\,}
\Per(\bar{\varphi}) \neq \emptyset$: Isometrically identify
an open disk $D \subset \Sigma$ with ${\mathcal S} =
\{ (x_0,\ldots,x_n)\, :\, x_0^2 + \cdots + x_n^2 = 1, x_n >
-1/2\}$. Let ${\bf g}$ be a $C^{\infty}$ metric on $\Sigma$ that
equals the round metric on $D \sim {\mathcal S}$. The geodesic flow of
${\bf g}$ possesses an open, invariant subset in $S_{\mathcal S}
\Sigma$ that consists entirely of periodic orbits, and these
orbits remain closed on the universal cover. One is led to ask

\medskip

\noindent
{\bf Question C}: {\em if $(\Sigma, {\bf g})$ is a compact
real-analytic riemannian manifold and $\pi_1(\Sigma)$ is anabelian, is
${\rm Int}\, \Lim(\bar{\phi}) = \emptyset$?}

\medskip

A connected component $B_i$ of $B$ is {\em plentiful} if the
inclusion $L_i \hookrightarrow S\Sigma$ has a finite-index image on
$\pi_1$. Let $(\bar{f},\bar{L},\bar{B})$ be the pullback of $(f,L,B)$
to $S\tilde{\Sigma}$.

\begin{thm} \label{thm:4.1}
Let $(\Sigma,{\bf g})$ be a compact $C^2$ finsler $n$-manifold whose
geodesic flow is $n$-semisimple with respect to $(f,L,B)$. Assume that
$\pi_1(\Sigma)$ is anabelian and $\Sigma$ is aspherical. If $B_i$ is
plentiful, then $\pi_k(\bar{B}_i)$ is non-trivial for some $k \geq 1$.
\end{thm}

\subsection{Background and Motivation} This paper is motivated by a question
posed by Kozlov~\cite{Kozlov:1979a}: which compact surfaces admit a
riemannian metric with an integrable geodesic flow?  Kozlov showed
that if the geodesic flow is analytic and has an additional analytic
first integral, then the surface's genus is at most one. Bolotin
subsequenty generalized Kozlov's argument to non-compact surfaces,
with an additional hypothesis on the behaviour of the metric at
infinity~\cite{Bolotin:1984a,Kozlov:1996a}. Further work by Bolotin
and Bolotin-Negrini shows that if $(\Sigma,{\bf g})$ is a compact
real-analytic surface of genus greater than $1$, then in a
neighbourhood of any non-trivial minimal periodic orbit of the
geodesic flow on $S\Sigma$ there is a
horseshoe~\cite{Bolotin:1993a,BolotinNegrini:1997a}. Ta\u{\i}manov~\cite{Taimanov:1988a,Taimanov:1995a}
generalized Kozlov's argument to higher dimensions, and obtained three
necessary conditions for a compact real-analytic manifold $\Sigma$ to
admit a real-analytically integrable geodesic flow: (1)
$\pi_1(\Sigma)$ must be almost abelian; (2) $\Sigma$'s first Betti
number $b$ is at most $\dim \Sigma$; and (3) there is an injection
of algebras $H^*(\T{b};\Q{}) \hookrightarrow
H^*(\Sigma;\Q{})$. Ta\u{\i}manov introduced what he called a {\em
geometrically simple} geodesic flow to prove these results; the
methods of the current paper are indebted to Ta\u{\i}manov's
conceptions.

Subsequently Paternain~\cite{Pat:1992a,Pat:1993a,Pat:1994a} introduced
the {\em entropy approach} to study integrable geodesic flows. He
showed that if a $C^{\infty}$ geodesic flow on a compact manifold
$\Sigma$ is integrable with first integrals that either (i) satisfy
Ito's non-degeneracy condition~\cite{Ito:1991a}; or (ii) generate a
$\T{n-1}$ action; or (iii) admit action-angle variables with
singularities, then the topological entropy of the geodesic flow must
vanish. By a result of Dinaburg and Bowen, $\pi_1(\Sigma)$ must be of
subexponential word growth. These results led Paternain to conjecture
that if $\Sigma$ admits a smoothly integrable geodesic flow (not
necessarily satisfying any of the hypotheses (i--iii) above), then
$\pi_1(\Sigma)$ is of polynomial word growth. Hence, $\pi_1(\Sigma)$
is almost nilpotent~\cite{Gromov:1981a}.

Note that Paternain's conclusions should be improvable. According to
Desolneux-Moulis~\cite{DesolneuxMoulis:1990a}, an Ito-nondegenerate
first-integral map admits a Whitney stratification. Thus, in cases (i)
and (iii), the critical-point set of the first-integral map is a
tamely-embedded polyhedron. If the image of the first-integral map
also admits a Whitney stratification -- and it almost certainly does
-- then the Kozlov-Ta\u{\i}manov theorem could be applied to conclude
that $\pi_1(\Sigma)$ is almost abelian.

The Kozlov-Ta\u{\i}manov theorem appeared to give a reasonably
complete characterization of manifolds with real-analytically
integrable geodesic flows, so the
examples~\cite{Butler:1999a,BolsinovTaimanov:2000a} were surprising
(see
also~\cite{Butler:2000a,Butler:2000b,Butler:2001a,BolsinovTaimanov:2000b}). In
essence, the current author showed integrability, in $C^{\infty}$
integrals, of the geodesic flows on $3$-dimensional manifolds with
$Nil$geometry. Bolsinov and Ta\u{\i}manov extended this construction
to a $3$-manifold with $Sol$geometry. The first example showed that
the necessary conditions for real-analytic integrablity derived by
Kozlov and Ta\u{\i}manov are not necessary for smooth integrability
(even of a real-analytic geodesic flow); the second examples showed
that integrable geodesic flows may have positive topological entropy
(even a real-analytic geodesic flow), thereby frustrating any simple
generalization of Paternain's work.  These examples were of elemental
importance in constructing the definitions of semisimplicity and
integrability offered in the first paragraph of the present paper.

\subsection{Outline}

Section~\ref{sec:kt} proves several technical lemmas based on the
technical hypotheses of Definition~\ref{def:fss}, then presents a
proof of the main lemma in~\cite{Taimanov:1988a}, suitably generalized
to the present setting. In addition to generalizing Ta\u{\i}manov's
lemma to non-compact manifolds, the proof shows the central importance
of condition FI2 of Definition~\ref{def:fss}. Section~\ref{sec:3mflds}
proves Theorem~\ref{thm:1}, Section~\ref{sec:neg} proves
Theorem~\ref{thm:3}. In Section~\ref{sec:ex}, the Butler and
Bolsinov-Ta\u{\i}manov examples are shown to be completely integrable
and semisimple. This section also shows that the geodesic flows on
associated infra-nil- and infra-solv-manifolds are completely
integrable and semisimple which suffices to prove Theorem~\ref{thm:2}.

\section{Definitions and Notation}

\noindent
Here are several useful conventions
\begin{enumerate}
\item[ $\bullet$ ] $\amalg$ denotes disjoint union;
\item[ $\bullet$ ] $\pi : \E \to \Sigma$ is the footpoint projection;
\item[ $\bullet$ ] $\tilde{\Sigma}$ (resp. $\hat{\Sigma}$) is the universal
cover (resp. a cover) of $\Sigma$;
\item[ $\bullet$ ] if $L \subset K$ then $\iota_L = \iota_{L,K}: L \to K$ is the inclusion map;
\item[ $\bullet$ ] a {\em geodesic flow} means the geodesic flow of a
complete $C^2$ finsler metric.
\end{enumerate}

Theorems~\ref{thm:1},\ref{thm:3},\ref{thm:4} are proven by applying a
more general result about {\em Hopf-Rinow} flows. As these results
have some independent applicability, it seems desirable to expose
them.

\subsection{HR flows}
Let $\pi : \E \to \Sigma$ be a smooth ($C^1$) fibre bundle with
compact fibres and let $\phi_t : \E \to \E$ be a flow. $\E$ may have a
boundary but $\Sigma$ is boundaryless and $\phi_t$ may only be a
local flow.

\begin{defn} \label{def:hr}
Let $q \in \Sigma$ and assume that for each non-trivial $[c] \in
\pi_1(\Sigma;q)$ there is a $p \in \pi^{-1}(q)$ and a $T > 0$ such that
$\gamma(t) := \pi \phi_{tT}(p)$, $0 \leq t \leq 1$, is a closed curve
homotopic to $c$; then $\phi_t$ is {\em Hopf-Rinow} over $q$. An
HR-flow is a flow that is Hopf-Rinow over $q$ for some $q \in \Sigma$.
\end{defn}

The curve $\gamma(t)$ will be called a {\em geodesic}. An HR-vector
field is one whose flow is HR. If two flows are orbitally equivalent
and one is HR, then so is the other. By the Hopf-Rinow theorem, the
geodesic flow of a complete $C^2$ finsler metric is an HR-flow on the
unit-sphere bundle~\cite{GiaquintaHildebrandt:1996a}. A skew product
over an HR-flow is also HR. A second class of Hopf-Rinow flows are
obtained as follows: let $\Sigma \subset M$ be an open submanifold
with a geodesically convex boundary and $\E = S_{\Sigma}M$. The
restriction of the geodesic flow to $\E$ defines a local flow that is
HR.

\subsection{$\FF$-semisimplicity}
Let us turn to a definition that abstracts some essential features of
complete integrability. A continuous surjection $f : L \to B$ is a
{\em fibration} if $f$ has the path-lifting property and the fibres of
$f$ are path-connected. If $B$ is paracompact and connected, then the
fibres of $f$ are of the same homotopy type, and so a ``typical''
fibre will be denoted by $F$~\cite{Whitehead:1978}.

\begin{defn} \label{def:fss} 
Let $\FF = \{f_i : L_i \stackrel{F_i}{\longrightarrow} B_i
\}_{i\ \in A}$ be a collection of fibrations. Let $\phi_t : \E
\to \E$ be an HR-flow, $L = \amalg_{i \in A}\ L_i$ and suppose that:
\begin{enumerate}
\item[(FI1)] $\Gamma = \E - L$ is closed, $\phi_t$-invariant and nowhere dense;
\item[(FI2)] for each $v \in \E$ and open neighbourhood $U \ni v$,
there is an open subset $W$, $v \in W \subseteq U$, such that $L \cap
W$ has finitely many path-connected components;
\item[(FI3)]  for each $i \in A$, $L_i$ is an open path-connected
component of $L$ and either $f_i \circ \phi_t = f_i$ or the inclusion
$F_i \hookrightarrow L_i$ induces an isomorphism on $\pi_1$.
\end{enumerate}
Then we will say that $\phi_t$ is {\em \fss}.
\end{defn}

If $\pi_1(F_i)$ is abelian for each $i \in A$, then $\phi_t$ will be
said to be abelian-\fss. Note that Definition~\ref{def:fss} does not
require the fibres $F_i$ to be compact.

\subsection{Related Definitions of Integrability} \label{ssec:connection}
Let $(P,\{,\})$ be a Poisson manifold. If $\F \subset C^2(P)$, let
$d\F_p = {\rm span}\, \{df_p\ :\ f \in \F\}$ and $Z(\F) = \{ f \in \F\
:\
\{\F,f\} \equiv 0\}$. When $\F$ is a Lie subalgebra of $C^{\infty}(P)$,
$Z(\F)$ is the centre of $\F$. Let $l = \sup
\dim d\F_p$, $k = \sup \dim dZ(\F)_p$. Say that a point $p \in P$ is
{\em strongly regular} for ${\bf f}$ if there is a saturated
neighbourhood $U \ni p$ and ${\bf f}|U$ is a trivial fibre-bundle
map.%
\footnote{If ${\bf f}$ is proper, then strong regularity is equivalent
to regularity.}
 A point $p \in P$ is $\F$-{\em regular} if there exists
$f_1,\ldots,f_l \in \F$ such that $p$ is strongly regular for the map
${\bf f} = (f_1,\ldots,f_l)$ and $f_1,\ldots,f_k \in Z(\F)$; if $p$ is
not $\F$-regular then it is $\F$-{\em critical}. Let $L(\F)$ be the
set of $\F$-regular points.

\begin{defn} \label{def:int}
$\F \subset C^2(P)$ is {\em tamely integrable} if
\begin{enumerate}
\item[ I1. ] $k + l = \dim P$;
\item[ I2. ] $L(\F)$ is an open and dense subset of
$P$;
\item[ I3. ] $P - L(\F)$ is a tamely-embedded polyhedron.
\end{enumerate}
A hamiltonian flow $\phi_t$ is tamely-$\F$-integrable if it enjoys a
tamely integrable set of first integrals.
\end{defn}

\noindent
See~\cite{BolsinovJovanovic:2001a} for an analogous definition and
further discussion. The conventional definitions of complete and
non-commutative integrability fit within the framework of
Definition~\ref{def:int}.

Let $(\Sigma,{\bf g})$ be a complete $C^2$ finslerian manifold. The
tangent bundle less its zero section, $\hat{T}\Sigma$, enjoys a
Poisson structure and the geodesic flow, $\phi_t$, is a $C^1$
hamiltonian flow on $\hat{T}\Sigma$ with $C^2$ hamiltonian $H$ (see
section 3.2 in~\cite{GiaquintaHildebrandt:1996a} for further
explanation).

\begin{thm} \label{thm:5}
Assume $\Sigma$ is compact. Then $\phi_t$ is tamely-$\F$-integrable
iff $\phi_t$ is integrable and semisimple.
\end{thm}

\begin{pf} If $\phi_t$ is tamely-$\F$-integrable, let $G$ denote the
abelian group of $C^1$ diffeomorphisms of $\hat{T} \Sigma$ generated
by the complete flows of $Y_h$, $h \in Z(\F)$. By the Sussman-Stefan
orbit theorem~\cite{Jurdjevic:1997a} and I1, the orbits of $G$ in
$L(\F)$ are embedded $C^1$ submanifolds. I1 and the properness of $H$
imply that for each $G$-orbit in $L(\F)$, there is a $G$-invariant
open neighbourhood, $U$, and an action of $\T{k}$ on $U$, such that
the $\T{k}$-orbits and $G$-orbits coincide. Thus each connected
component of $L(\F)$ is fibred by $\T{k}$-orbits. Therefore, there is
a $C^1$ atlas of $L(\F)$, ${\mathcal A} = \{ \varphi = (\theta,I) : U
\to \T{k} \times \R{l} \}$ which satisfies the
universal property that for all $\theta \in \T{k}$, $I \in \R{l}$ and
$1$-parameter subgroups $g^t$ of $G$: $\varphi \circ g^t
\circ \varphi^{-1}(\theta,I) = (\theta + t \xi(I), I)$ and $\xi : \R{l} \to
\R{k}$ is $C^1$. 

Let $L_o=L(\F)$, $B_o = L_o/G$ and $f_o : L_o \to B_o$ be the orbit
map. Since $\phi_t$ is a 1-parameter subgroup of $G$, it is integrable
with respect to $(f_o,L_o,B_o)$ . By I3 of Definition~\ref{def:int}
the complement of $L_o$ is a tamely-embedded polyhedron, so $\phi_t$
is semisimple with respect to $(f_o,L_o,B_o)$.

The opposite implication is straightforward. \qed
\end{pf}

Let $L = L_o \cap S\Sigma$, $B = f(L_o)$, $f = f_o | L$ and $\Gamma =
\Gamma_o \cap S\Sigma$. By hypothesis, $\Gamma_o = \hat{T}\Sigma - L_o$ is a
tamely-embedded polyhedron. Once can see that there is a triangulation
of $T\Sigma$ such that $\Gamma_o$ and $S\Sigma$ are
subcomplexes. Therefore $\phi_t | S\Sigma$ is integrable and
semisimple with respect to $(f,L,B)$. By defining $A$ to be the set of
connected components of $B$, and $\FF =
\{ f : L \to B \}$ it is immediately apparent that

\begin{thm} \label{thm:6}
If $\phi_t$ is tamely-$\F$-integrable, then $\phi_t$ is abelian-\fss.
\end{thm}

\section{An Extension of the Kozlov-Taimanov Theorem} \label{sec:kt}

This section generalizes the Kozlov-Ta\u{\i}manov theorem to \fss\ HR
flows. We start with a couple elementary lemmas. Let $d_H$ denote the
Hausdorff distance between two compacts sets.

\begin{lem} \label{lem:1} 
Assume that $\E = \Gamma \amalg L$ and that $L$ is dense and
satisfies (FI2). If $K \subset \E$ is compact, then given $\epsilon
> 0$, there is a bounded open set $U$, $K \subseteq U$,
$d_H(\overline{U},K)
\leq \epsilon$ such that $U \cap L$ has finitely many path-connected
components.
\end{lem} 

\begin{pf} 
Let $\epsilon > 0$ be given and let ${\mathcal O}$ denote the set
of open subsets of $\E$ whose diameter is less than $\epsilon/2$
and which intersect $L$ in finitely many path-connected components. By
(FI2), ${\mathcal O}$ is a covering of $\E$, hence of $K$. By the
compactness of $K$ there is a finite subcovering of $K$, which we will
denote by $K_1,\ldots,K_k \in {\mathcal O}$. We may assume that $K_i
\cap K \neq
\emptyset$ for all $i$. Then $U := \cup_{i=1}^k\ K_i$ is a
bounded open set containing $K$, and any point in $U$ is at most
$\epsilon$ away from a point in $K$. Since $L \cap K_i$ has finitely
many path-connected components, $U \cap L$ has finitely many
path-connected components. \qed
\end{pf}

\begin{lem}[Lifting Lemma] \label{lem:1b} 
Assume that $\phi_t : \E \to \E$ is a Hopf-Rinow flow that
satisfies (FI1 -- FI3). Let $p : \hat{\Sigma} \to \Sigma$ be a
covering of $\Sigma$. Then the flow $\hat{\phi}_t :
\hat{\E} \to \hat{\E}$ is Hopf-Rinow and satisfies (FI1 -- FI3).
\end{lem} 

\begin{pf} 
Clearly, the Hopf-Rinow property is satisfied when we pass to a
covering, so $\hat{\phi}_t$ is Hopf-Rinow. Let $\hat{\E} = p^* \E$ be
the pull-back of $\E$, and let $P : \hat{\E} \to \E$ denote the
covering map. We let $\hat{\Gamma} = P^{-1}(\Gamma)$ and $\hat{L} =
P^{-1}(L)$. Since $P$ is a local homeomorphism and conditions (FI1)
and (FI2) are purely local, they are obviously satisfied. Let $C$
denote the connected components of the set $\hat{L}$, so that $\hat{L}
= \amalg_{j \in C}\, \hat{L}_j$ (we abuse notation and let $\hat{L}_j
= j$ for each $j \in C$). Clearly, for each $j \in C$ there is an $i
\in A$ such that $P(\hat{L}_j) = L_i$; we define $\check{f}_j :
\hat{L}_j \to B_i$ by $\check{f}_j = f_i \circ P$.  Since $f_i$
has the path-lifting property and $P$ is a local homeomorphism,
$\check{f}_j$ has the path-lifting property~\cite{Whitehead:1978}; and
clearly $\check{f}_j$ is surjective.

Let $\hat{\pi}_j : \hat{B}_j \to B_i$ be a covering space of $B_i$
such that $\im \hat{\pi}_{j,*} = \im \check{f}_{j,*}$. By the usual
properties of covering spaces, there is a continuous surjective map
$\hat{f}_j$ such that
$$ 
\begin{array}{lcl}
\hat{L}_j & \stackrel{\hat{f}_j}{\longrightarrow} & \hat{B}_j\\
          & \stackrel{\check{f}_j}{\searrow}      & \downarrow \hat{\pi}_j\\
          &                                       & B_i
\end{array}
$$
commutes. Since $\check{f}_j$ has the path-lifting property and
$\hat{\pi}_j$ is a local homemorphism, $\hat{f}_j$ has the
path-lifting property. By construction, the fibres of $\hat{f}_j$ are
path-connected, and since $\hat{f}_j$ is continuous and surjective, it
is a fibration.

Let $\hat{F}_j = \hat{f}_j^{-1}(b)$ for some $b \in \hat{B}_j$.  It
follows that $\hat{\FF} = \{ \hat{f}_j : \hat{L}_j
\stackrel{\hat{F}_j}{\longrightarrow} \hat{B}_j \}_{j \in C}$ is a
collection of fibrations that satisfies (FI1--FI3). \qed
\end{pf}

\begin{lem}[c.f. Theorem 1~\cite{Taimanov:1988a}] \label{lem:2}
Assume that $\phi_t : \E \to \E$ is a Hopf-Rinow flow that satisfies
{\rm (FI1)} and {\rm (FI2)}. Then there is an $i \in A$ such that the
map $$\pi_1(L_i) \to \pi_1(\E) \to \pi_1(\Sigma)$$ has a finite-index
image.
\end{lem} 

\begin{pf}
Assume that the Hopf-Rinow flow $\phi_t : \E
\to \E$ satisfies (FI1) and (FI2) of Definition~\ref{def:fss}. Let
$L := \cup_{\alpha \in A}\, L_{\alpha}$. Let $q \in \Sigma$ be a point
at which $\phi_t$ satisfies the Hopf-Rinow property. Let $\E_q =
\pi^{-1}(q)$ and let $Q$ and $P$ be open disks containing $q$ such
that $\overline{Q} \subset P$. Let $U = \pi^{-1}(Q)$, $W =
\pi^{-1}(P)$. The contractibility of $Q$ (resp. $P$) means that $U$
(resp. $W$) is topologically trivial, so $V \simeq Q \times \E_q$
(resp. $W \simeq P \times \E_q$). Since $\E_q$ is compact,
lemma~\ref{lem:1} implies that there is a bounded open set $V$ such
that $\overline{U} \subset V \subset W$ and $L \cap V$ has a finite
number of path-connected components. Let $K_1,\ldots,K_{\omega}$ be an
enumeration of these path-connected components.

Because $K_i \subset W$, $\pi(K_i) \subset P$ so $\pi(K_i)$ can be
contracted to the point $q$ in $W$ for all $i$. In addition, for each
$i = 1,\ldots,\omega$, $K_i$ is a path-connected subset of $L$. By
(FI2) each path-connected subset of $L$ lies in a unique component
$L_{\alpha}$. That is, for each $i \in
\{1,\ldots,\omega\}$ there is a unique $\alpha_i \in A$ such that $K_i
\subset L_{\alpha_i}$. Finally, by (FI1) $L$ is dense in
$\E$ so $\cup_{i=1}^{\omega}\ K_i$ is a dense subset of $V$.

For each $j = 1,\ldots,\omega$, select a point $v_j \in K_j$ and let
$a_j : [0,1] \to P$ be a continuous curve that joins $q$ to $q_j :=
\pi(v_j)$. Let $a_j^*(t) := a_j(1-t)$ be the curve traversed in the
opposite sense. By construction, each $q_j$ lies in the contractible
set $P$.

\noindent
Let us agree to call a curve $\gamma(t) = \pi \circ \phi_t(x)$, for $0
\leq t \leq T$ a {\em geodesic}; $\pi(x)$ is the footpoint of the
geodesic $\gamma$.

Because $\phi_t$ is a Hopf-Rinow flow, for each non-trivial homotopy
class $[\tilde{c}] \in \pi_1(\Sigma;q)$ there is a geodesic in
$[\tilde{c}]$ with footpoint $q$. Let $u \in \E_q$ be the initial
condition of such a geodesic $\gamma$, and let its length be $T$. The
initial condition $u$ may lie in the singular set $\Gamma$, but by the
continuity in initial conditions of $\phi_t$, for each $\epsilon > 0$
there exists a $\delta > 0$ such that if $v \in \E$ and $d(u,v) <
\delta$, then $d( \phi_t(u), \phi_t(v) ) < \epsilon$ for all $t \in [0,T]$. By the
openness of $V$, the density of $L \cap V$ in $V$, and the invariance
of $L$ (FI2), for all $\epsilon > 0$ sufficiently small there is a $v
\in L \cap V$ such that $\phi_T(v) \in L \cap V$ and so $\pi \circ
\phi_T(v) \in P$. Because $L \cap V = \cup_{i=1}^{\omega}\ K_i$, there are $i,j \in \{1,\ldots,\omega\}$ such that
$v \in K_i$ and $\phi_T(v) \in K_j$. Indeed, by the
$\phi_t$-invariance of each $L_{\alpha}$ (FI2), there exists an
$\alpha$ such that $K_i,K_j \subset L_{\alpha}$.

Let now $C$ be the arc in $\E$ that consists of an arc $s : [0,1]
\to K_i$ joining $v_i$ to $v$, followed by the arc obtained by
following the trajectory $\phi_t(v)$ for $0 \leq t \leq T$, followed
by an arc $e : [0,1] \to K_j$ joining $\phi_T(v)$ to $v_j$. Let $c$ be
the arc in $\Sigma$ obtained by concatenating $a_i$, $\pi \circ C$ and
$a_j^*$. The contractibility of $P$ implies that the curve $c$ is
homotopic to the geodesic arc through $u$, namely $\gamma(t) = \pi
\circ \phi_t(u)$ for $t \in [0,T]$. Therefore $c \in [\tilde{c}]$.

Take the collection of all such arcs $C$ in $L$ constructed in the
previous paragraph. These arcs generate a groupoid: if $C$ ends in
$(K_i,v_i)$ and $D$ begins in $(K_i,v_i)$, then their product
(concatenation) $D*C$ is defined. Modulo homotopies in $L$ that fix
end points, the operation $*$ is associative, and so the equivalence
classes of arcs $C$ generates a groupoid $\G$. The collection of
homotopy classes of arcs $C$ that begin and end in $(K_i,v_i)$
generates a group $\G_i$ for each $i=1,\ldots,\omega$. $\G_i$ is a
subgroup of $\pi_1(L_{\alpha_i}; v_i)$, where $K_i \subset
L_{\alpha_i}$.

Observe that: (i) if $C_{ji}$ is an arc beginning in $(K_i,v_i)$ and
ending in $(K_j,v_j)$ then $C_{ji}^{-1} \G_j C_{ji} = \G_i$; (ii)
there is a subset $\{C_{ij}\} \subset \G$ of cardinality no greater
than $\omega^2$ such that for any element $C \in \G$ there is a
$C_{ij}$ such that $C_{ij} C \in \G_i$ for some $i,j$. Claim (i)
follows from the observation that if $[C] \in \G_j$, then
$[C]*[C_{ji}]$ is a homotopy class (relative to endpoints) of loops
that begin at $v_j$, end at $v_i$ and remain in $L$. Therefore,
$[C_{ji}]^{-1} * [C] * [C_{ji}]$ is a homotopy class (relative to
endpoints) of loops that begin at $v_i$, end at $v_i$ and remain in
$L$, i.e. $[C_{ji}]^{-1} * [C] * [C_{ji}] \in \G_i$. Claim (ii)
follows from the finiteness of the number of components $K_i$.

The map $C \in \G$ $\to$ $[c] = [a_j^* (\pi \circ C) a_i] \in
\pi_1(\Sigma;q)$ is an epimorphism $s$ of groupoids induced by
the maps $L \stackrel{\iota_L}{\longrightarrow} \E
\stackrel{\pi}{\longrightarrow} \Sigma$. The restriction $s | \G_i$ is
a group homomorphism. By (ii),
$$\pi_1(\Sigma;q) = \cup_{1\leq i,j\leq \omega}\ c_{ij} H_i$$ where
$c_{ij} = s(C_{ij})$ and $H_i = s(\G_i)$. Therefore $\pi_1(\Sigma;q)$
is a finite union of cosets of subgroups. To prove that at least one
of the subgroups $H_i$ is of finite index in $\pi_1(\Sigma;q)$ it
remains to observe

\begin{lem}[\cite{Taimanov:1988a}]
Suppose that a group $G = \cup_{1\leq i\leq \alpha,\ 1 \leq j
\leq \beta}\ c_{ij} H_i$ where $c_{1,1},\ldots,$ $c_{\alpha,\beta} \in
G$ and $H_1,\ldots,H_{\alpha}$ are subgroups of $G$. Then there is at
least one subgroup $H_i$ of finite index in $G$.
\end{lem}

\noindent
Since $H_i \leq (\pi \iota_{L_{\alpha_i}})_*\,
\pi_1(L_{\alpha_i};v_i)$, this completes the proof. \qed
\end{pf}

Combining the lifting lemma with Lemma~\ref{lem:2} shows that,
possibly after passing to a finite covering, one can assume that for
some $i$, $\pi_1(L_i) \to \pi_1(\Sigma)$ is epimorphic.

\section{$3$- and $4$-manifolds} \label{sec:3mflds}

Recall that if an integrable flow $\phi_t : M \to M$ is semisimple
with respect to $\FF = (f,L,B)$, we say that $\FF' = (f',L',B')$ is a
{\em refinement} of $\FF$ if $B'$ is open and dense in $B$, $L' =
f^{-1}(B')$, $f' = f|L'$ and $\Gamma' = M - L'$ is a nowhere-dense
tamely-embedded polyhedron.

\def\trac{{tractable}}

\begin{defn} \label{def:niceref}
$\FF = (f,L,B)$ is {\em \trac} if, for each connected component $f_i :
L_i \to B_i$, either
\begin{enumerate}
\item[1.] there is a compact codimension-1 submanifold $N_i \subset B_i$ such
that the inclusion of $W_i = f_i^{-1}(N_i) \hookrightarrow L_i$ is
epimorphic on $\pi_1$; or
\item[2.] there is a component $L_j$ satisfying {\rm 1.} and a map $r : L_i
\to L_j$ such that $\iota_{L_i,M}$ is homotopic to $\iota_{L_j,M}
\circ r$.
\end{enumerate}
\end{defn}

\noindent
{\em Remark}: Condition 2. implies that $\im \iota_{L_i,M*} \subset
\im \iota_{L_j,M*}$. Hence, from the point of view of the fundamental
group, one only need concern oneself with those components $L_j$ that
satisfy Condition 1.

\begin{lem} \label{lem:niceref}
If $\phi_t : M \to M$ is semisimple with respect to $\FF = (f,L,B)$
and $M$ is compact, then there is a {\trac} refinement $\FF' =
(f',L',B')$ of $\FF$.
\end{lem}

\def\G{{\bf G}}
\def\U{{\bf U}}
\def\K{{\bf K}}
\begin{pf} We break the proof into 3 cases.

\noindent
{\em Case 1}: $\dim B = 1$ is trivial.

\noindent
{\em Case 2}: $\dim B = 2$.

\noindent
{\em Claim}: $B$ is homeomorphic to a compact surface $S$ punctured at
a finite number of points.

\noindent
To prove the claim, it suffices to show that $B$ has finitely-many
ends. If $B$ has no ends, it is compact and the claim is trivial. If
$B$ is a one-ended surface whose boundary is $\T{1}$, then $B$ is
homeomorphic to $\T{1} \times [0,1)$. If $B$ has finitely-many ends,
then $B = C \cup E$ where $C$ is a compact surface with boundary and
$E$ is a regular neighbourhood of infinity. Hence each component of
$E$ is a one-ended surface with a $\T{1}$ boundary, i.e. each
component is homeomorphic to $\T{1} \times [0,1)$. This proves the
claim if $B$ has finitely-many ends.

Since $f : L \to B$ is proper, it induces a bijection of ends. Let $K$
be a polyhedral neighbourhood of $\Gamma$, and let $K_n$ be a regular
neighbourhood of $\Gamma$ in the $n$-th barycentric subdivision of
$K$. Let $\Gamma$ have $\kappa$ components. Compactness of $M$ implies
that $\kappa$ is finite. Hence $C_n = K_n \cap L$ has $\kappa$
components for all $n$ sufficiently large. Number these components
$C^j_n$ such that $C^j_{n+1} \subset C^j_n$ for all $j$ and all $n$
sufficiently large. This shows that $L$, hence $B$, has $\kappa <
\infty$ ends. This establishes the claim.

Let $S$ be as in the claim. Triangulate $S$ so that each puncture
point is a barycentre of a 2-simplex. This triangulation induces a
decomposition of $B$ into a 1-skeleton $B^{(1)}$ and a finite union of
open 2-disks and open cylinders. Let $B' = B - B^{(1)}$, $L' =
f^{-1}(B')$ and $f' = f|L'$. For each component $B'_i$ of $B'$, let
$N_i$ be an embedded circle such that $\pi_1(N_i) \to \pi_1(B'_i)$ is
epimorphic. Clearly $W_i = f^{-1}(N_i) \hookrightarrow L'_i$ is
surjective on $\pi_1$, so $(f',L',B')$ is tractable.

From~\cite{Putz}, the proper fibration $f : L \to B$ is
triangulable. Hence, it can be assumed that $f : L \to B$ is a
simplicial map of PL manifolds. From this, it is clear that it can be
assumed that $f^{-1}(B^{(1)})$ is a compact subcomplex of $L$ and
hence a tamely-embedded polyhedral subset of $M$. Since $\Gamma$ and
$f^{-1}(B^{(1)})$ are tamely-embedded polyhedra which are separated by
an open set, $\Gamma' = \Gamma \cup f^{-1}(B^{(1)})$ is a
tamely-embedded polyhedron. Hence $(f',L',B')$ is a tractable
refinement of $(f,L,B)$.

\noindent
{\em Case 3}: $\dim B \geq 3$.

\noindent
{\em Step 1}: Let $\G \subset \K$ be simplicial complexes with
polytopes $|\G| \subset |\K|$ such that $\G$ is a full subcomplex in
$\K$'s codimension-1 skeleton. Let $\U$ be the subcomplex of $\K$
obtained by deleting all simplices in $\K$ with a vertex in $\G$. From
the proof of Lefschetz duality, there is a deformation retraction of
$|\K| - |\G|$ onto $|\U|$~\cite{Munkres:1984}.

\def\CC{{\mathcal C}}

\noindent
{\em Step 2}: Let $K$ be a polyhedral neighbourhood of $\Gamma$ in
$M$. From step 1, $L$ admits a deformation retraction onto $L_o = L -
{\rm Int}\, K$. Let $\Phi : L \times I \to L$ be this deformation
retraction.

From~\cite{Putz}, the proper fibration $f : L \to B$ is
triangulable. Hence, it can be assumed that $f : L \to B$ is a
simplicial map of PL manifolds and that $\Phi$ is a simplicial
map. Let $B_o = f(L_o)$, $B_1$ be a regular neighbourhood of $B_o$ and
$B_+ = {\rm Int}\, B_1$. Let $B_- = B - B_1$, so that $B$ is the
disjoint union of $B_+$, $B_-$ and $\partial B_1$. Let $L_{\pm} =
f^{-1}(B_{\pm})$.

Let $B' = B_- \cup B_+$, $L' = f^{-1}(B')$ and $f' = f|L'$. Since
$\Gamma$ and $f^{-1}(\partial B_1)$ are tamely-embedded polyhedra that
are separated by open sets, $\Gamma'$ is a tamely-embedded
polyhedron. Hence $(f',L',B')$ is a refinement of $(f,L,B)$.

For the next two paragraphs, assume that $B$ and $B_-$ are
connected. Then $B_{+}$ is connected.

Let $B^{(1)}_+ \subset B_+$ be the union of 1-simplices of $B_1$ that
lie entirely in $B_+$. Let $N_+ \subset B_+$ be the boundary of a
regular neighbourhood of $B^{(1)}_+$ and let $W_+ =
f^{-1}(N_+)$. Since $\dim B \geq 3$, $N_+$ is connected and the
inclusion $N_+
\hookrightarrow B_+$ is epimorphic on $\pi_1$; hence $W_+
\hookrightarrow L_+$ is epimorphic on $\pi_1$. This proves $L_+$
satisfies 1. in Definition~\ref{def:niceref}.

On the other hand, let $r : L_- \to L_+$ be $r = \iota_{L_o,L_+} \,
\circ \,
\Phi_1 \, \circ \, \iota_{L_-,L}$ where $\Phi_t(\bullet) =
\Phi(\bullet,t)$. Since $\Phi_0 = id_L$,
$\iota_{L_+,L} \circ r$ is homotopic in $L$ to $\iota_{L_-,L}$. Hence
$\iota_{L_+,M} \circ r$ is homotopic in $M$ to $\iota_{L_-,M}$.  This
proves $L_-$ satisfies 2. in Definition~\ref{def:niceref}.

If $B$ or $B_-$ is not connected, the previous two paragraphs can
be applied componentwise. This proves the Lemma. \qed
\end{pf}

\begin{pf*}{Proof of Theorems~\ref{thm:1} and~\ref{thm:3}.} 
Let $(\Sigma,{\bf g})$ be a compact finslerian $n$-manifold and assume
that the geodesic flow $\phi_t : S\Sigma \to S\Sigma$ is integrable
and semisimple with respect to $(f,L,B)$. The dimension of $B$
(resp. $L$) is $l$ (resp. $k+l$) and $k+l=2n-1$.

By Lemma~\ref{lem:niceref}, it can be assumed that $\FF = (f,L,B)$ is
tractable. By Lemma~\ref{lem:2}, there is a connected component $L_i$
of $L$ such that $\im(\pi \iota_{L_i})_*$ is of finite index in
$\pi_1(\Sigma)$. By the remark following Definition~\ref{def:niceref},
it can be assumed that $L_i$ satisfies Condition 1., hence there is a
compact codimension-1 manifold $W_i \subset L_i$ that fibres over a
compact codimension-1 manifold $N_i \subset B_i$ such that
$\iota_{W_i,L_i*}$ is epimorphic on $\pi_1$. Let us suppress the index
$i$ in the following discussion.

The homotopy exact sequence $\pi_2(N)
\stackrel{ \partial_*}{\to} \pi_1(\T{k}) \to \pi_1(W) \to \pi_1(N)
\to 1$, implies there is a short exact sequence $$1 \to \Omega \to \pi_1(W)
\to \pi_1(N) \to 1,$$ where $\Omega$ is abelian of
rank at most $k$. Thus $\pi_1(W)$ is polycyclic if $\pi_1(N)$ is
polycyclic. The step length of $\pi_1(W)$ is at most $k$ plus the step
length of $\pi_1(N)$. Since $\pi_1(L) = \im
\iota_{W,L*}$, these comments apply to $\pi_1(L)$, too.

If $\phi_t$ is $n$-semisimple, then $k=n$ and $l=n-1$ so $\dim N =
n-2$. When $n=3$, $\pi_1(N) \simeq \Z{}$, so $\pi_1(L)$ is polycyclic
of step length at most $4$. When $n=4$, $N$ is a compact surface. From
the hypothesis of Theorem~\ref{thm:3}, $N$ covers the projective
plane or the Klein bottle, so $\pi_1(N)$ contains a finite-index copy
of $\Z{s}$, $s \leq 2$. Thus, $\pi_1(L)$ contains a finite-index
polycyclic subgroup of step length at most $6$. \qed
\end{pf*}

\section{Anabelian fundamental groups and $\FF$-semisimplicity} \label{sec:neg}

\begin{lem}  \label{lem:nm1}
Let $\pi_1(\Sigma)$ be anabelian, $L \subseteq \E$, and $f : L
\stackrel{F}{\longrightarrow} B$ be a fibration. If $\pi_1(F)$ is
abelian, then either:
\begin{enumerate}
\item[1.] the image of $L \stackrel{\pi \iota_L}{\longrightarrow}
\Sigma $ on $\pi_1$ is not of
finite index; or
\item[2.] the composite map $F \stackrel{\iota_{F}}{\hookrightarrow} L
\stackrel{\pi \iota_L}{\longrightarrow} \Sigma $ is trivial on $\pi_1$
and the homomorphism $(\pi \iota_{L})_* :
\pi_1(L) \to \pi_1(\Sigma)$ factors through a homomorphism
$\xi : \pi_1(B) \to \pi_1(\Sigma)$.
\end{enumerate}
\end{lem}

\begin{pf} Because $f : L \to B$ is a fibration, there is a
pair of intersecting horizontal and vertical exact sequences
$$
\begin{array}{cccccccc}
          &
& \ \ \ker (\pi \iota_L)_*\\
          &
& \downarrow \ \ \ \ \ \ \ \\                      
\pi_1(F)  & \stackrel{\iota_{F,*}}{\longrightarrow} 
& \pi_1(L) \ \ \ \    & \stackrel{f_*}{\longrightarrow} & \pi_1(B)\\
          &                                         & \ \ \downarrow (\pi \iota_L)_* 
&         &  \downarrow id \\
          &                                        &
 \pi_1(\Sigma) \ \ \ \ & \stackrel{?\xi}{\longleftarrow} & \pi_1(B),
\end{array} 
$$
and $?\xi$ indicates that $\xi$ remains to be defined. If $\im (\pi
\iota_L)_*$ is of finite index, then $\im (\pi
\iota_L \iota_F)_*$ is an almost normal abelian subgroup of
$\pi_1(\Sigma)$. Hence it is trivial, so $\im \iota_{F,*} \subset \ker
(\pi \iota_L)_*$ which suffices to define $\xi$ as the diagram
suggests. \qed
\end{pf}

\begin{pf*}{Proof of Theorem~\ref{thm:4}.}
From the definition of integrability, each point on a
$\phi_t$-invariant torus is recurrent. For each $v
\in L$, let $U$ be an open neighbourhood of $v$ that is contractible
in $L$.  Because $v \in L$ is recurrent, there is a sequence $t_k
\to \infty$ as $k \to \infty$ such that $\phi_{t_{k}}(v), \to v$
as $k \to \infty$. Without loss of generality, it may be assumed that
$\phi_{t_k}(v) \in U$ for all $k$. Let $\gamma_k : \T{1} \to \E$ be
the closed loop obtained by concatenating the orbit segment
$\phi_t(v)$ for $0 \leq t \leq t_k$ with an arc $E_k$ from
$\phi_{t_k}(v)$ to $v$. It is clear that $E_k$ can be chosen to lie in
$U$ and have length($E_k$) $\to 0$ as $k \to \infty$. Since $U$ is
contractible in $L$ and $f$ is $\phi_t$-invariant, $f
\circ \gamma_k$ is null-homotopic in $B$.

Let $i \in A$ be such that $\im(\pi \iota_{L_i})_*$ is of finite index
in $\pi_1(\Sigma)$. By Lemma~\ref{lem:nm1}, $(\pi \iota_{L_i})_*$
factors through a homomorphism $\xi_i : \pi_1(B_i) \to
\pi_1(\Sigma)$. Then, the loop $\gamma_k$ is freely homotopic to a loop in
$F_i$, so $\gamma_k$ is null-homotopic in $\E$. Therefore, any lift of
$\gamma_k$ to $\bar{\E}$ is also closed. But this implies that any
$\bar{v} \in \bar{L}_i$ in the fibre over $L_i$ is recurrent for
$\bar{\phi}_t$. Since $\bar{L}_i$ is open, $\Lim(\bar{\phi})$ has a
non-empty interior.

For the second part of the theorem, we use the method of toroidal
surgeries~\cite{Bogo:1996a}. Let $\Phi=(\theta,I) : U \to
\T{k} \times \R{l}$ be a $C^1$ diffeomorphism with proper inverse that
conjugates $\phi_t$ with the translation-type flow $T_t(\theta,I) =
(\theta + t \xi(I), I)$. Let $D \subset \R{l}$ be a small open ball
contained in the image of $I$ and let $\Xi : D \to {\bf P}^{k-1}$ be
the frequency map $\Xi(I) = [\xi_1(I) : \cdots : \xi_k(I) ]$. Since $k
\leq l+1$, in any $C^1$-neighbourhood of $\Xi$, there is an $\Omega =
[\omega_1(I): \cdots : \omega_k(I)]$ such that $\Omega$ is a
submersion on $D$ and the support of $\Xi - \Omega$ is the closure of
$D$. Let $S_t(\theta,I) = (\theta+t\omega(I),I)$ and define
$$
\eta_t(P) = \left\{ \begin{array}{cc}
\phi_t(P)             & {\rm if\ } P \not\in U,\\
\Phi^{-1} S_t \Phi(P) & {\rm if\ } P \in U.
\end{array}
\right.
$$
It is clear that $S_t$, hence $\eta_t$, is $C^1$. It is also clear
that $\Per(\eta)$ contains an open set. From the first part of the
theorem, it follows that $\Per(\bar{\eta})$ contains an open set.

Finally, $\eta_t$ is orbitally equivalent to the geodesic flow of a
finsler metric by the arguments in sections 3.2--3.3
of~\cite{Butler:2001a}. \qed
\end{pf*}

\begin{pf*}{Proof of Theorem~\ref{thm:4.1}.} 
Let $n = \dim \Sigma$. Since $\phi_t$ is completely integrable $\dim B
= n-1$.

If $B_i$ is plentiful, the lifting lemma implies that one may assume,
possibly after passing to a finite covering of $\Sigma$, that
$\im(\pi\iota_{L_i})_* = \pi_1(\Sigma)$ and that $\Sigma$ is
orientable. Let us suppress the subscript $i$ in the remainder of the
proof.

Assume that $\pi_k(\bar{B}) = 1$ for all $k$. Then $B$ is aspherical
and the epimorphism $\xi : \pi_1(B) \to \pi_1(\Sigma)$ of
Lemma~\ref{lem:nm1} is an isomorphism. Since $B$ and $\Sigma$ are
aspherical and their fundamental groups are isomorphic, they are
weakly homotopy equivalent, hence their singular homology groups are
isomorphic (Corollary V.4.6 in~\cite{Whitehead:1978}). But $\Sigma$ is
a compact orientable manifold, so its top homology group $H_n(\Sigma)$
is non-zero, while $H_n(B) = 0$. Absurd. \qed
\end{pf*}

\section{Examples} \label{sec:ex}

\subsection{The Poisson Geometry of $T^* G$}
The proof of Theorem~\ref{thm:2} exploits an underlying Lie-theoretic
structure. The principal machinery used to construct integrals will be
the {\em momentum map}.

\def\g{{\mathfrak g}}

Let $G$ be a real Lie group with Lie algebra $\g$. The dual space,
$\g^*$, has a Poisson bracket defined as follows: The derivative of a
smooth function on $\g^*$ at a point is naturally identified with an
element in $\g$. The Poisson bracket on $\g^*$ is defined by
\begin{equation} \label{eq:ps}
\{f,h\}_{\g^*}(\mu) := -\langle \mu , [ df_{\mu}, dh_{\mu}] \rangle,
\end{equation}
for all $f,h \in C^{\infty}(\g^*)$ and $\mu \in \g^*$. A function $f
\in C^{\infty}(\g^*)$ is a {\em Casimir} if $E_f = \{.,f\}$ is trivial. The
set of Casimirs is precisely $Z(C^{\infty}(\g^*))$. The vector field
$E_f$ is called the {\em Euler vector field}.

The map
\begin{equation}
\psi(g,\mu) := \Ad{g}^*\, \mu
\end{equation}
is called the momentum map of $G$'s left-action on $\g^*$. A
right-invariant vector field on $G$ is of the form $\xi^R_{G}(g) = d_e
R_g\, \xi$, for some $\xi \in \g$, and all $g \in
\G$. The cotangent lift of $\xi^R_{G}$ has the hamiltonian function
$h_{\xi}(g,\mu) = \langle \xi^R_{\G}(g), \mu \rangle$, which equals
$\langle \psi(g,\mu), \xi \rangle$. It is known that

\noindent
{\em The map $\psi : T^* G \to \g^*$ is the momentum map of $G$'s
left-action on $T^* G$. The map $\omega(g,\mu) = \mu$ is the momentum
map of $G$'s right action on $T^* G$. Both maps are submersions.}

The canonical Poisson structure on $T^* G$, $\{,\}_{T^* G}$, is
related to that on $\g$, $\{,\}_{\g^*}$ as follows: $\{,\}_{T^* \G}$
is right (resp. left) invariant, so the Poisson bracket of right
(resp. left) invariant functions is again right (resp. left)
invariant. If ${\mathcal R}$ (resp. ${\mathcal L}$) denotes the right
(resp. left) invariant functions smooth functions on $T^* G$, then
${\mathcal R} = \psi^* C^{\infty}(\g^*)$ (resp. ${\mathcal L} =
\omega^* C^{\infty}(\g^*)$), and $\psi^*$ (resp. $\omega^*$) is a Lie
algebra isomorphism (resp. anti-isomorphism). In addition, because
right and left multiplication commute these two subalgebras commute:
$\{ {\mathcal R}, {\mathcal L} \}_{T^* G} \equiv 0$. Because $\psi$ is
a Poisson map, we will abuse notation and use $\{,\}$ to denote either
$\{,\}_{T^* G}$ or $\{,\}_{\g}$, depending on the context.

Therefore, if $H \in {\mathcal L}$, then $H$ Poisson commutes with all
hamiltonians $F \in {\mathcal R}$ so $F$ is a first integral of the
hamiltonian vector field $Y_H$ on $T^* G$. In addition, the projection
map $\omega : T^* \G \to \g^*$ also satisfies $\omega_* Y_H =
-E_H$. Thus, if $k \in C^{\infty}(\g^*)$ is a first integral of $E_H$,
then $k \circ \omega$ is a first integral of $Y_H$.

Consequently, to prove integrability of $Y_H$ on $T^* G$, it is useful
to: (1) find sufficiently many functions in ${\mathcal R}$; and (2)
find integrals of $X_H$ on $\g^*$. With luck, the sum of these two
subalgebras of integrals will be sufficient for integrability. In the
event that we wish to study $Y_H$ on $T^* (E \backslash G)$, we need
to find sufficiently many functions in ${\mathcal R}^{E}$.

\subsection{An integrability theorem}

\newcommand \isom [1] {{\rm Isom(}#1 {\rm )}}
\def\A{{\mathcal A}}
\def\B{{\mathcal B}}
\def\a{{\mathfrak a}}
\def\b{{\mathfrak b}}
\def\EE{{\mathcal E}}
\def\AA{{\mathcal A}}
\def\BB{{\mathcal B}}
\def\C{{\bf C}}

\noindent
{\em Standing Hypothesis}: It will be assumed throughout that the
Casimirs of $\g^*$ separate $G$'s coadjoing orbits.

For a left-invariant metric ${\bf g}$ on $G$, let $\isom{{\bf g}}$ be
the isometry group of ${\bf g}$ and let $O({\bf g})$ be the group of
automorphisms of $\g$ that are also ${\bf g}$-orthogonal. Let $E$ be a
discrete, torsion-free subgroup of $\isom{{\bf g}}$ that acts freely
and uniformly discretely on $G$; $\Sigma = E \backslash G$ is the
quotient manifold and ${\bf h}$ is the metric on $\Sigma$ induced by
${\bf g}$. There is an exact sequence
$$1 \to N \to E \to F \to 1,$$ where $N = E \cap G$ and $F = E / N$ is
isomorphic to a finite subgroup of $O({\bf g})$. 

For $\alpha \in \Aut{G}$, let $\beta = (d_e \alpha')^{-1}$ be the
induced linear isomorphism on $\g^*$ and let $T^* \alpha$ be the
symplectomorphism of $T^* G$ induced by $\alpha$. A calculation shows
that $\psi \circ T^* \alpha = \beta \circ \psi.$ The natural
right-action of $G$ on $T^* G$ extends to an action of $G
\star \Aut{G}$ on $T^* G$, and this action factors through to an
action on $\g^*$. Let $\EE$ be the subgroup of GL($\g^*$) induced by
$E < G \star \Aut{G}$.

\begin{thm} \label{thm:gint}
Assume that $\g^* = \g^*_r \amalg \g^*_s$, $\BB
\subset C^{\omega}(\g^*)$ and $\AA \subset C^{\omega}(\g^*_r)$ satisfy
\begin{enumerate}
\item[ F1. ] $\g^*_s$ is a nowhere dense, analytic, $\EE$-invariant
set;
\item[ F2. ] $\BB$ is an integrable subalgebra of $C^{\omega}(\g^*)^F$ containing the
hamiltonian of ${\bf k}$;
\item[ F3. ] $\AA$ is an integrable subalgebra of
$C^{\omega}(\g^*_r)^{\EE}$.
\end{enumerate}

Then: if $\Sigma$ is compact, then the geodesic flow of ${\bf k}$ is
integrable and semisimple.
\end{thm}

\begin{pf} Let $\tilde{\b} = \omega^* \BB$ and $\tilde{\a} = \psi^*
\AA$. Since $\tilde{\b}$ is left-invariant and $F$-invariant, it induces a subalgebra
$\b$ on $T^*\Sigma$. Similarly, $\tilde{\a}$ is $\EE$-invariant by F3,
so it induces a subalgebra $\a$ on $T^* \Sigma$. Without changing F2
and F3, it can be assumed that both $\AA$ and $\BB$ contain the
Casimirs of $\g^*_r$. Let $c$ be the index of $G$, and $n$ be the
dimension of $G$, so that the dimension of a coadjoint orbit in
$\g^*_r$ is $n-c$ and $G$ enjoys $c$ independent Casimirs on each
coadjoint orbit in $\g^*_r$, by the standing hypothesis. Since $F$ is
finite, it can be assumed that each of these Casimirs is also
$F$-invariant.

\def\Reg{{\rm Reg}}

Let $T^*\Sigma_r = \Sigma \times_{F} \g^*_r$. Due to analyticity and
the hypothesis that both $\AA$ and $\BB$ are integrable subalgebras,
there exists $\lambda : T^* \Sigma_r \to \R{c+b}$ (resp. $\rho : T^*
\Sigma_r \to \R{c+a}$) whose regular level set $\Reg(\lambda)$ (resp. 
$\Reg(\rho)$) is an everywhere-dense analytic set such that:
\begin{enumerate}
\item[P1.] rank $d\lambda_P = b+c$ for $P \in \Reg(\lambda)$ (resp. rank $d\rho_P
= a+c$ for $P \in \Reg(\rho)$);
\item[P2.] $\lambda_j = \rho_j$ for $j=1,\ldots,c$;
\item[P3.] $\lambda_j$ for $j=1,\ldots,c$ are induced by $F$-invariant Casimirs of $\g^*$;
\item[P4.] $\{ \lambda_i, \lambda_j \} = 0$ for $i=1,\ldots,b+c$ and
$j=1,\ldots,b_o+c$ (resp. $\{ \rho_i, \rho_j \} = 0$ for
$i=1,\ldots,a+c$ and $j=1,\ldots,a_o+c$);
\item[P5.] $a_o+a=b_o+b=n-c$.
\end{enumerate}
Let ${\bf f}_j = \lambda_j$ for $j=1,\ldots,c+b_o$ and ${\bf f}_{j+c+b_o}
= \rho_{j+c}$ for $j=1,\ldots,a_o$, ${\bf f}_{j+a_o+c} =
\lambda_{j}$ for $j=b_o+1,\ldots,b$ and ${\bf f}_{j+b+c} = \rho_{j+c}$
for $j=a_o+1,\ldots,a$. Since $\tilde{\a} \cap \tilde{\b}$ is a set of
bi-invariant functions on $T^* G$, P1-P3 imply that $\Reg({\bf f}) =
\Reg(\lambda) \cap \Reg(\rho)$ is an everywhere-dense analytic subset
of $T^* \Sigma_r$ hence of $T^* \Sigma$.

P4-P5 imply that $\{ {\bf f}_i, {\bf f}_j \} = 0$ for
$i=1,\ldots,a_o+b_o+c$ and $j=1,\ldots,a+b+c$. Let $\F$ be the
subalgebra of $C^{\infty}(T^*
\Sigma)$ generated by the ${\bf f}$-pullbacks of compactly-supported
smooth functions whose support is contained in the interior of $\im
{\bf f}$. Then $\F$ satisfies:
\begin{enumerate}
\item[ J1. ] $k = \dim dZ(\F)_P = a_o+b_o+c$ for all $P \in \Reg({\bf f})$;
\item[ J2. ] $l = \dim d\F_P = a+b+c$ for all $P \in \Reg({\bf f})$;
\item[ J3. ] $\Gamma = T^* \Sigma - \Reg({\bf f})$ is a closed, nowhere-dense
analytic set.
\end{enumerate}

\noindent
From P5, $k+l = 2n$ so J1-J2 imply that I1 of Definition~\ref{def:int}
is satisfied. Since $\Gamma$ is an analytic subset of $T^* \Sigma$,
there is a triangulation of $(T^*
\Sigma,\Gamma)$~\cite{Lojasiewicz:1964a}.  Hence $\Gamma$ is a
tamely-embedded polyhedron, so $\F$ is a tamely-integrable
algebra. Since $\F$ is a set of integrals of the geodesic flow of
${\bf k}$, Theorem~\ref{thm:5} implies that the geodesic flow is
integrable and semisimple. \qed
\end{pf}

\subsubsection{Construction of the integrable subalgebra $\AA$}

The following two sets of conditions imply the existence of an
integrable subalgebra $\AA \subset C^{\omega}(\g^*_r)^{\EE}$. Note that
H1-H3 is a special case of G1-G4. The proofs are straightforward and
left to the reader. $C^{\infty}_o$ is the set of smooth functions with
compact support.

Assume that $\g^* = \g^*_r \amalg \g^*_s$ is the disjoint union of two
sets such that either
\begin{enumerate}
\item[ H1. ] $\g^*_s$ is a closed, nowhere-dense, $\EE$-invariant
real-analytic set;
\item[ H2. ] there is a real-analytic fibration $\C :
\g^*_r \to B$ such that $\Ad{G}^*$ acts transitively on the fibres of
$\C$;
\item[ H3. ] for each $b \in B$, $\Ad{N}^*$ acts freely and uniformly
discretely on $\C^{-1}(b)$;
\end{enumerate}

or
\begin{enumerate}
\item[ G1. ] $\g^*_s$ is a closed, nowhere-dense, $\EE$-invariant
real-analytic set;
\item[ G2. ] there is a real-analytic $G$-manifold $V$ and
$G$-equivariant real-analytic fibrations $\g^*_r
\stackrel{{\bf p}}{\to} V \stackrel{\C}{\to} B$;
\item[ G3. ] there is a normal subgroup $N_{stab}$ of $N$ such that for each $b
\in B$, $N/N_{stab}$ acts freely and uniformly discretely on
$\C^{-1}(b)$;
\item[ G4. ] $\dim V + \dim d Z( {\bf p}^* C_o^{\infty}(V) ) = n + c$.
\end{enumerate}

\section{Applications}

This section proves Theorem~\ref{thm:2}. The geometric $3$-manifolds
${\bf E}^3$, $Nil$ and $Sol$ will play an important role.

\noindent
\subsection{${\bf E}^3$} 
Let $G$ be the Lie group $\R{3}$ and let ${\bf g}$ be a left-invariant
metric on $G$; ${\bf E}^3 = (G, {\bf g})$. The isometry group of ${\bf
E}^3$ is naturally isomorphic to $G \star O(3)$. The Bieberbach
Theorem says that if $E$ is a uniformly discrete, torsion-free
cocompact subgroup of ${\rm Isom}({\bf E}^3)$, then there is an exact
sequence
$$1 \to N \to E \to F \to 1,$$ where $N = G \cap E$ is the maximal
abelian subgroup of $E$ -- which is isomorphic to $\Z{3}$ -- and $F$
is a finite subgroup of $O(3)$. Relative to the obvious trivialization
of $T^* G$, the momentum map of $G$'s right action on $T^* G$ is
$$\psi(h,p) = \Ad{h}^* p = p.$$ Let $\g^*_r$ be the subset of $\g^*$
on which $F$ acts freely, let $B = \g^*_r$ and let ${\bf C} =
id$. Conditions H1--H3 are obviously satisfied. Since $G$ is abelian,
its coadjoint orbits are all points, so the standing hypothesis is
trivially satisfied. In this case, $a_o = b_o = 0$ and $c = 3$. Thus,

\begin{thm} \label{thm:e3}
Let $\Sigma = E \backslash G$ where $E$ is a uniformly discrete,
torsion-free and cocompact subgroup of isometries of ${\bf g}$. The
geodesic flow of the metric induced by ${\bf g}$ on $\Sigma$ is
completely integrable and semisimple.
\end{thm}

Note that the orbifold $\g^* / F$ is naturally coordinatized by the
$F$-invariant polynomials on $\g^*$. The geodesic flows in Theorem
\ref{thm:e3} are integrable with real-analytic first
integrals. The monodromy of the Liouville foliation is also naturally
isomorphic to $F$.

The analogous theorem for the geometric $3$-manifolds that are modeled
on $S^2 \times {\bf E}^{1}$ and $S^3$ is proven in an identical manner
and will be omitted.

\noindent
\subsection{$Nil$} c.f.~\cite{Butler:1999a,Butler:2000a,Butler:2000b}.  
Let $G = Nil$ denote the set $\R{3}$ with the multiplication rule:
$(x,y,z)*(x',y',z') = (x+x', y+y', z+z'+\frac{1}{2}(xy'-x'y))$. $G$ is
the $3$-dimensional Heisenberg group whose center is $Z(G) = \{
(0,0,z) \}$. The $1$-forms $\alpha = dx$, $\beta = dy$ and $\gamma =
dz - \frac{1}{2}( xdy - ydx )$ are left-invariant and form a basis of
$\g^*$. The dual basis ${\mathcal X}, {\mathcal Y}$ and ${\mathcal Z}$
of $\g$ satisfies $\exp( x {\mathcal X} + y {\mathcal Y} + z{\mathcal
Z}) = (x,y,z)$ for all $x,y,z \in \R{}$. Thus, $\g$ and $G$ can be
identified with $\R{3}$ in the obvious way, and in this coordinate
system the exponential map is just the identity. Let
\begin{equation}
{\bf g} = \sum_{\sigma \in \{\alpha,\beta,\gamma\}}
\sigma \otimes \sigma.
\end{equation}

The isometry group of ${\bf g}$ is the semi-direct product of $G$ with
the subgroup of automorphisms whose derivative at $id$ is
orthogonal. Relative to the $(x,y,z)$ coordinate system the
automorphism group acts as linear transformations and
\[ O( {\bf g} ) = 
\left\{ \left[ \begin{array}{ccc} \cos \theta & -\sin \theta & 0\\  
s\sin \theta & s\cos \theta & 0\\ 0 & 0 & s
\end{array} \right]\, :\, s = \pm 1, \theta \in \R{} \right\}.\] 
$O({\bf g})$ is a maximal compact subgroup of $\Aut{G}$. It is known
that if $E$ is a uniformly discrete, torsion-free subgroup of $G \star
\Aut{G}$, then $E$ is conjugate to a subgroup of $G \star O({\bf g})$
and the maximal normal nilpotent subgroup $N$ of $E$ is simultaneously
conjugated to a subgroup of the group generated by the elements
$(1,0,0)$, $(0,2,0)$ and $(0,0,1)$~\cite{Dekimpe:1996a}. So it may be
assumed that $E$ is a subgroup of $G \star O({\bf g})$ with these
properties. Then, there is a commutative diagram
\[ \left.
\begin{array}{lcccccccr}
           &     & 1          & \to & 1          &     &            &     & \\
           &     & \downarrow &     & \downarrow &     &            &     & \\
1          & \to & Z          & \to & Z          & \to & 1          &     & \\
           &     & \downarrow &     & \downarrow &     & \downarrow &     & \\
1          & \to & N          & \to & E          & \to & F          & \to & 1\\
           &     & \downarrow &     & \downarrow &     & \downarrow &     &  \\
1          & \to & L          & \to & Q          & \to & F_o        & \to & 1\\
           &     & \downarrow &     & \downarrow &     & \downarrow &     &  \\
           &     & 1          &     & 1          &     & 1.         &     &
\end{array}
\right.
\] where $Q$ is a crystallographic group of motions of the plane, $L$
is a lattice subgroup of $\R{2}$, $F_o$ is a finite group of linear
isometries of the plane, and $Z = N \cap Z(G)$ is the center of $N$
and $F$ is isomorphic to a finite subgroup of $O({\bf g})$. All of the
maps in the diagram are the obvious ones induced from the exact
sequence $1 \to Z(G) \to G \to G/Z(G) \to 1$ where $G/Z(G)$ is
identified with $\R{2}$.

Write a covector $p \in T^*_h G$ as $p = p_{\alpha} \alpha + p_{\beta}
\beta + p_{\gamma} \gamma$; this amounts to trivializing $T^* G$ with
respect to the left action of $G$. The momentum map of $G$'s right
action on $T^* G$ is then:
\[\psi(h,p) = \Ad{h}^* p = (p_{\alpha} +
\frac{1}{2}y p_{\gamma}) \alpha + (p_{\beta} - \frac{1}{2} x
p_{\gamma}) \beta + p_{\gamma} \gamma,\] where $h = (x,y,z)$. It is
clear that the coadjoint orbits of $G$ on $\g^*$ are planes $\{
p_{\gamma} = c \}$ for non-zero $c$ and single points for $c = 0$. Let
$\g^*_r = \{ p \in \g^*\ :\ p_{\gamma} \neq 0 \}$ and let $\g^*_s$ be
the complement of $\g^*_r$. It is clear from the explicit description
of $O({\bf g})$ and its action on $\g$, that $\g^*_r$ is $O({\bf
g})$-invariant. Let $B = \R{ \times }$ and define ${\bf C} : \g^*_r \to B$ by
\begin{equation} \label{eq:c}
{\bf C}(p) := p_{\gamma}.
\end{equation} 
Then H1 and H2 are satisfied.

Since $N/Z$ is isomorphic to a lattice subgroup $L$ of the plane, it
is clear from the explicit description of $\Ad{h}^* p$, that
$\Ad{N}^*$ acts freely and uniformly discretely on each fibre of ${\bf
C}$. Hence H3 is satisfied.

Let the quadratic form on $\g^*$ induced by ${\bf g}$ be denoted by
$g$. It is given by
\begin{equation}
g(p) = \frac{1}{2} \sum_{\sigma \in \{\alpha,\beta,\gamma\}}
p_{\sigma}^2.
\end{equation}
Let $\B = {\rm span}\, \{ g, p_{\gamma} \}$. Since the dimension of
$G$ is $3$ and its index is $1$, $\B$ is an integrable subalgebra of
$C^{\omega}(\g^*)$. Thus

\begin{thm} \label{thm:nil}
Let $\Sigma = E \backslash G$ where $E$ is a uniformly discrete,
torsion-free cocompact subgroup of isometries of ${\bf g}$. The
geodesic flow of the metric induced by ${\bf g}$ on $\Sigma$ is
integrable and semisimple.
\end{thm}

\noindent
The following discussion proves the $3$-semisimplicity of the geodesic flow
of ${\bf g}$.

\subsubsection{Chern classes and monodromy} The choice of $\C$ in
Equation~\ref{eq:c} is not unique and it turns out that alternative
choices of $\C$ have interesting geometric properties. Let's explain
with the simplest case case where $E = N$. 
\begin{enumerate}
\item[ Case 1. ] ${\bf C} = p_{\gamma}$.\\
The map
\begin{equation}
\zeta(N h,p) = (\ g(p),\ p_{\gamma},\
-\frac{2 p_{\beta}}{p_{\gamma}} + x\ {\rm mod}\ 1,\frac{2
 p_{\alpha}}{p_{\gamma}} + y\ {\rm mod}\ 1\, )
\end{equation}
is a real-analytic mapping defined of $T^* \Sigma_r \to \R{+} \times
(\R{ \times } ) \times \T{1} \times \T{1}$. The first two components of $\zeta$
Poisson commute with all components, and $\zeta$ is a proper
submersion except on the set $\Sigma \times
\{p_{\alpha}=p_{\beta}=0\}$. Thus, the geodesic flow of ${\bf g}$ is non-commutatively
integrable. It is clear that it is also semisimple. As $\zeta$ is
derived from the canonical first-integral map $(Nh,p) \to ( g(p),
\Ad{N}^* \psi(h,p) )$ from $T^* \Sigma_r \to \R{+} \times \left(
\g^*_r / \Ad{N}^* \right)$, this essentially reproves Theorem~\ref{thm:nil}.

Let $L = \Sigma \times \{ p_{\alpha}^2+p_{\beta}^2 > 0 \}$ be the
regular-point set of $\zeta$. $L$ has the structure of a $\T{2}$-fibre
bundle, so we can compute its monodromy and Chern class. Let $L_+ :=
\Sigma \times \{ p_{\alpha}^2+p_{\beta}^2 > 0, p_{\gamma} > 0 \}$ be
one of the two connected components of $L$ and let $\zeta | L_+$ be
denoted by $f : L_+ \to B$ where $B = \R{+} \times (\R{ \times })
\times \T{2}$. $L_+$ retracts onto $L_0 := \Sigma \times \{
p_{\alpha}^2+p_{\beta}^2 = 1, p_{\gamma} = 1 \} \simeq \Sigma \times
\T{1}$. Since $\Sigma$ is a principal $\T{1}$-bundle over $\T{2}$ with
projection map $\pi(N h) = Z(G) N h$ -- i.e. $\pi(N(x,y,z)) = (x\ {\rm
mod}\ 1,y\ {\rm mod}\ 1)$ -- $f | L_0$ equals $(N h,\theta\ {\rm
mod}\, 1) \to (2,1,\sigma(\theta) + \pi(N h))$, where $\sigma(\theta)
= (-2\sin
\theta, 2 \cos \theta)$. Since $\sigma$ is
null-homotopic, $f|L_0$ is homotopic to $(N h, \theta\ {\rm mod}\, 2\pi)
\to \pi(N h)$. Thus $f|L_0$ is homotopic to the composition of the
canonical projections $\Sigma \times \T{1}
\to \Sigma \to \T{2}$. Therefore, the monodromy group of $f : L_+ \to
B$ is trivial and the Chern class of $f$ is naturally identified with
the Chern class of $\pi \times id$, which is non-trivial.

\item[ Case 2. ] ${\bf C} = p_{\gamma} \gamma \oplus \left(
\frac{2p_{\alpha}}{p_{\gamma}} + y + \Z{} \right)\alpha$.\\
In this case,
\begin{equation}
\xi(N h,p) = (\ g(p),\ p_{\gamma},\
-\frac{2 p_{\beta}}{p_{\gamma}} + x\ {\rm mod}\ 1 )
\end{equation}
is real-analytic surjection of $T^* \Sigma_r \to \R{+} \times (\R{
\times } ) \times \T{1}$. The components of $\xi$ Poisson commute and
$\xi$ is a proper submersion except on the set $\Sigma \times
\{p_{\alpha}=p_{\beta}=0\}$. Thus the geodesic flow of ${\bf g}$ is completely
integrable on $T^* \Sigma$. It is clear that it is also semisimple.

Let $L_+$ be as above, and let $\ell = \xi| L_+$. Clearly, $\ell$ is a
proper lagrangian fibration whose image is a manifold with trivial
second ${\rm \check{C}ech}$ cohomology group. Hence the Chern class of
$\ell$ is trivial. The subgroup ${\mathcal V} =
\{ (0,y,z) \}$ is normal in $G$, and this normal subgroup endows the
manifold $\Sigma$ with the structure of a $\T{2}$ bundle over $\T{1}$,
with projection map defined by $\Pi(Nh) = {\mathcal V}Nh$; that is
$\Pi(N(x,y,z)) = x\, {\rm mod}\, 1$. Arguments similar to those above
show that $\ell | L_0$ is homotopic to the composition of canonical
projections $\Sigma \times \T{1} \to \Sigma \to
\T{1}$. Thus, the monodromy of $\ell | L_0$ is equal to that of $\Pi
\times id$, which is non-trivial.

Thus, we have an example of an integrable system that is naturally
tangent to a lagrangian foliation which has non-trivial monodromy but
a trivial Chern class, and it is also tangent to an isotropic
foliation which has trivial monodromy but a non-trivial Chern class.
\end{enumerate}

\subsection{$Sol$} See~\cite{BolsinovTaimanov:2000a,BolsinovTaimanov:2000b}.
Let $G$ be $\R{3}$ equipped with the multiplication
$(x,y,z)*(x',y',z') = (x+x', \exp(x)y' + y, \exp(-x)z' + z)$. Let
${\mathcal V}$ be the normal subgroup $\{ (0,y,z) \}$; ${\mathcal V}$
is naturally isomorphic to $\R{2}$ and $G = \R{+} \star \R{2}$ is the
semi-direct product of $\R{+}$ with $\R{2}$ and is a 3-dimensional
solvable Lie group. The $1$-forms $\alpha = dx$, $\beta = \exp(-x)\,
dy$, and $\gamma = \exp(x)\, dz$ are left-invariant and span
$\g^*$. Let
\begin{equation}
{\bf g} = \sum_{\sigma \in \{\alpha,\beta,\gamma\}}
\sigma \otimes \sigma
\end{equation}
be a left-invariant riemannian metric on $G$. The isometry group of
$G$ is equal to $G \star O({\bf g})$ where $O({\bf g})$ is the group
of orthogonal automorphisms of $G$. The automorphism group of $G$
acts, in the $(x,y,z)$ coordinate system, as a group of linear
transformations. The group of orthogonal automorphisms is isomorphic
to the dihedral group of order eight with the elements
\[ O( {\bf g} ) = 
\left\langle \left[ \begin{array}{ccc}  1 &   0 & 0\\  
 0 & \pm 1 & 0\\ 0 & 0 & \pm 1
\end{array} \right] \ , \  \left[ \begin{array}{ccc} -1 & 0 & 0\\  
 0 & 0 & \pm 1\\ 0 & \pm 1 & 0
\end{array} \right] \right\rangle,\] where the sign of the diagonal elements are independent of each other~\cite{Scott:1983a}.

Let $d$ be a positive square-free integer and let $K = \Q{}(\sqrt{d})$
be the totally real quadratic number field obtained by adjoining
$\sqrt{d}$ to $\Q{}$. Let $\sigma$ be the unique non-trivial
automorphism of $K$ that fixes $\Q{}$ and maps $\sqrt{d}$ to
$-\sqrt{d}$. Let ${\mathcal O}_K$ be the ring of integers of $K$,
${\mathcal I}$ an ideal of ${\mathcal O}_K$ and let $u > 0$ be a
non-trivial unit of ${\mathcal O}_K$. The additive group of ${\mathcal
I}$ is a $u$ module, so let $\Delta = \langle u \rangle \star
{\mathcal I}$. Define an embedding of $\Delta$ into $G$ by
$$(u,a) \mapsto ( \log u, a, \sigma(a) ), $$ for all $a \in {\mathcal
I}$. It is elementary to verify that this is a group embedding and
that the image, $N$, is a discrete cocompact subgroup of $G$.

Write a covector $p \in T^*_h G$ as $p = p_{\alpha} \alpha + p_{\beta}
\beta + p_{\gamma} \gamma$; this amounts to trivializing $T^* G$ with
respect to the left action of $G$. The momentum map of $G$'s right
action on $T^* G$ is
\[\psi(h,p) = \Ad{h}^* p = (p_{\alpha} + y \exp( x )\, p_{\beta} + z
\exp( -x )\, p_{\gamma})
\alpha + \exp( x )\, p_{\beta} \beta + \exp( -x )\, p_{\gamma}
\gamma,\] where $h = (x,y,z)$. It is clear that the regular coadjoint orbits
of $G$ on $\g^*$ are connected components of the level sets of the
Casimir $\kappa = p_{\beta} p_{\gamma}$. Let $\g^*_r = \{ p \in
\g^*\ :\ \kappa(p) \neq 0 \}$ and let $\g^*_s$ be the complement of
$\g^*_r$. Since $O({\bf g})$ is a group of automorphisms of $G$,
$\g^*_r$ is $O({\bf g})$-invariant. Let $V = \R{ \times } \beta \oplus
\R{ \times } \gamma$ and let ${\bf p}(p) := p_{\beta} \beta +
p_{\gamma} \gamma$. The action of $\Ad{G}^*$ on $\g^*_r$ factors
through the map ${\bf p}$. Let $B = (\R{ \times }) \times \Z{}_2$ and
define ${\bf C} : V
\to B$ by
\[ {\bf C}(p_{\beta} \beta + p_{\gamma}\gamma) := (p_{\beta}p_{\gamma}, {\rm
sign}(p_{\beta}) ).\] Thus G1 and G2 are satisfied.

Let $N_{stab} = N \cap {\mathcal V}$. From the explicit description of
the coadjoint action, it is clear that $N/N_{stab} \simeq \langle u
\rangle$ acts freely and uniformly discretely on the fibres of $\C$. Thus G3 is true.

Note that $V$ is $2$-dimensional, ${\bf p}^* C_o^{\infty}(V)$ is an
abelian subalgebra of $C_o^{\infty}(\g^*)$ and $\dim G = 3$, ${\rm
ind}\, G = 1$, so G4 is satisfied. Let $\BB = {\rm span} \{ g, \kappa
\}$ where $g$ is the quadratic form on $\g^*$ induced by ${\bf
g}$. This shows that

\begin{thm} \label{thm:sol}
Let $\Sigma = E \backslash G$ where $E$ is a uniformly discrete,
torsion-free cocompact group of isometries of ${\bf g}$. The geodesic
flow of ${\bf g}$ is completely integrable and semisimple.
\end{thm}

\begin{pf*}{Proof of Theorem~\ref{thm:2}.} Let $\Sigma$ be a $3$-manifold
with $\pi_1(\Sigma)$ infinite polycyclic and $\pi_2(\Sigma) = 0$. From
Evans and Moser's theorem~\cite{EvansMoser:1972a}, $\pi_1(\Sigma)$ is
isomorphic to one of the following:
\begin{enumerate}
\item $\Z{}$, $\Z{2}$ or $\mathfrak K$, the fundamental group of the
Klein bottle;
\item an extension $1 \to A \to \pi_1(\Sigma) \to \Z{} \to 1$ where $A
= \Z{2}$ or $\mathfrak K$;
\item an amalgamation $\langle a,b,x,y\, |\, bab^{-1} = a^{-1},yxy{-1}
= x^{-1}, a=x^p y^{2q}, b^2 = x^r y^{2s} \rangle$ where $p,q,r,s$ are
integers and $|ps-rq| = 1$;
\item an extension of a group in (2), with $A=\Z{2}$, by a finite group of automorphisms.
\end{enumerate}

The groups in (1) are not the $\pi_1$ of a compact, boundaryless
$3$-manifold with $\pi_2 = 0$. For $\Z{2}$ (hence $\mathfrak K$), this
is Reidemeister's theorem. Similarly, if $\Sigma$ is a compact,
boundaryless $3$-manifold with $\pi_1(\Sigma) = \Z{}$ and
$\pi_2(\Sigma) = 0$, then $\Sigma$ is a $K(\Z{},1)$-space, hence
homotopy equivalent to $\T{1}$. But by Poincar\'e duality the second
Betti number of $\Sigma$ is $1$. Absurd.

In case (2), irreducibility of $\Sigma$ plus Theorem 3
of~\cite{HempelJaco} imply that $\Sigma$ fibres over $\T{1}$ with
fibre $\T{2}$ or the Klein bottle -- hence $\Sigma$ admits flat, $Nil$
or $Sol$ geometry by Theorem 5.3 of~\cite{Scott:1983a}. 

In cases (3) and (4), $\pi_1(\Sigma)$ contains a finite-index subgroup
of type (2). Hence $\Sigma$ is finitely covered by a $\T{2}$-bundle
over $\T{1}$, so $\Sigma$ admits flat, $Nil$ or $Sol$ geometry by
Theorem 5.3 of~\cite{Scott:1983a}.

Thus, by Theorems~\ref{thm:e3}--\ref{thm:sol}, $\Sigma$ admits a
real-analytic riemannian metric whose geodesic flow is
$3$-semisimple. \qed
\end{pf*}

\bibliographystyle{plain} 
\bibliography{bibliography} 
\end{document}